\documentclass{article}

\usepackage{amssymb}
\usepackage{amsmath}
\usepackage[dvips]{graphicx}

\begin{document}


\noindent \textbf{                                  }

\noindent \textbf{}

\noindent \textbf{}

\noindent \textbf{}

\noindent \textbf{}
\begin{center}
 \textbf{INVERSE INDEFINITE SPECTRAL PROBLEM FOR SECOND }

 \textbf{ORDER DIFFERENTIAL OPERATOR WITH }

 \textbf{COMPLEX PERIODOC COEFFICIENTS  }
 \end{center}
\noindent \textbf{}
\begin{center}
 \textbf{R.F. Efendiev }
\end{center}
\begin{center}
\textbf{Institute Applied Mathematics, Baku State University,}

\textbf{Z.Khalilov, 23, AZ1148, Baku, Azerbaijan, }
\end{center}
\begin{center}
rakibaz@yahoo.com\textbf{ \textit{  }}
\end{center}
\textbf{\textit{}}

\textbf{\textit{}}

\noindent{ABSTRACT.}

\noindent

     The inverse problem for the Sturm- Liouville operator with complex periodic
potential and discontinuous coefficients on the axis is studied.
Main characteristics of the fundamental solutions are
investigated, the spectrum of the operator is studied.  We give
formulation of the inverse problem, prove a uniqueness theorem and
provide a constructive procedure for the solution of the inverse
problem.


spectral problems; Continuous spectrum;


\noindent\textbf{{INTRODUCTION.}}

\textbf{\underbar{}}

\indent        We consider the differential equation

\begin{equation} \label{GrindEQ__1_}
-y''\left(x\right)+\, q\left(x\right)y\left(x\right)=\lambda ^{2} \rho \left(x\right)y\left(x\right)
\end{equation}
in the space $L_{2} \left(-\infty ,+\infty \right)$ where the prime denotes the derivative with respect to the space coordinate and assume that the potential  $q\left(x\right)$ is of the form

\begin{equation} \label{GrindEQ__2_}
q(x)=\sum _{n=1}^{\infty }q_{n} e^{inx}  ,
\end{equation}
the condition$\sum _{n=1}^{\infty }\left|q_{n} \right|^{2} =q<\infty  $ is satisfies, $\lambda $ is a complex number, and

\begin{equation} \label{GrindEQ__3_}
\rho \left(x\right)=\left\{\, \, \begin{array}{c} {1} \\ {-\beta ^{2} } \end{array}\right. \, \, \, \begin{array}{c} {for} \\ {for} \end{array}\, \, \, \, \begin{array}{c} {x\ge 0,} \\ {x<0.\, \, } \end{array}
\end{equation}
      The function $\rho \left(x\right)$ is called density function, and the function $q\left(x\right)$ is called the potential function of equation \eqref{GrindEQ__1_}. If$\rho \left(x\right)=1$, then equation \eqref{GrindEQ__1_} is called the potential equation.  The potential equation often is met in physical, technical and astronomical problems. Generalized Legendre  equation, degenerate   hypergeometrical equation, Bessel's equation and also Mathieu equation after suitable substitution coincide with potential equation \eqref{GrindEQ__1_}[1,p.374]-[2].

\noindent        As a rule, such problems are connected with discontinuous properties of materials. Inverse problems of spectral analysis consist of recovering operators from their spectral characteristics.

\noindent In this paper we will study the spectrum and also solve the inverse problem for singular non-self-adjoint operator.  As the coefficient allows bounded analytic continuation to the upper half-plane of the complex plane$z=x+it$, we can conduct detailed analysis of problem \eqref{GrindEQ__1_}-\eqref{GrindEQ__3_}.

\noindent        Our investigation was stimulated by M.G.Casymov's paper [3] where he first considered this potential, and his co-workers [4]. Especially, we would like to indicate the paper V. Guillemin, A. Uribe [5] where the potential $q(x)=\sum _{n=1}^{\infty }q_{n} e^{inx}  $ plays a vital part for solving the KdV equation. Later in 1990 the results obtained in [3] were extended by Pastur L.A., Tkachenko V.A [6].

\noindent The inverse problems connected with potentials of the form $q(x)=\sum _{n=1}^{\infty }q_{n} e^{inx}  $

\noindent   where$\sum _{n=1}^{\infty }\left|q_{n} \right|=q<\infty  $, were considered in [7-9].

\noindent The paper consists of three parts.

\noindent    In part 1 we study the properties of fundamental system of solutions of equation \eqref{GrindEQ__1_}. The spectrum of problem \eqref{GrindEQ__1_}-\eqref{GrindEQ__3_} is investigated in part 2. In part 3 we give a formulation of the inverse problem, prove a uniqueness theorem and provide a constructive procedure for the solution of the inverse problem.

\textbf{\underbar{1.  REPRESENTATION OF FUNDAMENTAL SOLUTIONS.}}

\noindent \textbf{\underbar{}}

Here we study the solutions of the main equation

\noindent

\[-y''\left(x\right)+\, q\left(x\right)y\left(x\right)=\lambda ^{2} \rho \left(x\right)y\left(x\right)\]
that will be convenient in future.

\noindent  We first consider the solutions $f_{1}^{+} \left(x,\lambda \right)$ and $f_{2}^{+} \left(x,\lambda \right)$, determined by the conditions at infinity

\[\, \, \mathop{\lim }\limits_{Imx\to \infty } f_{1}^{+} \left(x,\lambda \right)e^{-i\lambda x} =1,\]

\[\, \, \, \mathop{\lim }\limits_{Imx\to \infty } f_{2}^{+} \left(x,\lambda \right)e^{-\beta \lambda x} =1.\]
We can prove the existence of these solutions if the condition $\sum _{n=1}^{\infty }\left|q_{n} \right|^{2} =q<\infty  $ is fulfilled for the potential. This will be unique restriction on the potential and later on we'll consider it to be fulfilled.

\noindent \textbf{Theorem 1.  }Let \textit{$q(x)$} be of the form \eqref{GrindEQ__2_} and  \textit{$\rho \left(x\right)$ }satisfy condition \eqref{GrindEQ__3_}. Then equation \eqref{GrindEQ__1_} has special solutions of the form

\begin{equation} \label{GrindEQ__4_}
f_{1}^{+} (x,\lambda )=e^{i\lambda x} \left(1+\sum _{n=1}^{\infty }\frac{1}{n+2\lambda } \sum _{\alpha =n}^{\infty }V_{n\alpha } e^{i\alpha x}   \right)\, ,\, \, \, \, \, \, \, \, for\, \, x\ge 0,
\end{equation}

\begin{equation} \label{GrindEQ__5_}
f_{2}^{+} (x,\lambda )=e^{\lambda \beta x} \left(1+\sum _{n=1}^{\infty }\frac{1}{n-2i\lambda \beta } \sum _{\alpha =n}^{\infty }V_{n\alpha } e^{i\alpha x}   \right)\, ,\, \, \, \, \, \, for\, \, \, x<0.
\end{equation}
where the numbers $V_{n\alpha } $ are determined from the following recurrent relations

\begin{equation} \label{GrindEQ__6_}
\alpha (\alpha -n)V_{n\alpha } +\sum _{s=n}^{\alpha -1}q_{\alpha -s} V_{ns} =0 ,\, \, \, \, \, \, \, \, \, \, \, 1\le n<\alpha ,
\end{equation}

\begin{equation} \label{GrindEQ__7_}
\alpha \sum _{n=1}^{\alpha }V_{n\alpha } +q_{\alpha } =0 ,
\end{equation}
 and the series

\begin{equation} \label{GrindEQ__8_}
\sum _{n=1}^{\infty }\frac{1}{n} \sum _{\alpha =n}^{\infty }\alpha \left|V_{n\alpha } \right|
\end{equation}
converges.

\noindent The proof of the theorem is similar to the proof of [6] and therefore we don't cite it here.

\noindent \textbf{Remark1:  } If $\lambda \ne -\frac{n}{2} $ and$Im\lambda \ge 0$, then$f_{1}^{+} \left(x,\lambda \right)\in L_{2} \left(0,+\infty \right)$.

\noindent \textbf{Remark2:}  If $\lambda \ne -\frac{in}{2\beta } $ and$Re\lambda \ge 0$, then$f_{2}^{+} \left(x,\lambda \right)\in L_{2} \left(-\infty ,0\right)$.

\noindent  Extending $f_{1}^{+} \left(x,\lambda \right)$and $f_{2}^{+} \left(x,\lambda \right)$ as solutions of equation \eqref{GrindEQ__1_} on\textit{$x<0$} and\textit{ $x\ge 0$ }respectively and using the conjunction conditions

\begin{equation} \label{GrindEQ__9_}
\begin{array}{l} {y\left(0+\right)=y(0-),} \\ {y'\left(0+\right)=y'(0-),} \end{array}
\end{equation}
we can prove the following lemma.

\textbf{Lemma 1: } $f_{1}^{+} \left(x,\lambda \right)$ and $f_{2}^{+} \left(x,\lambda \right)$  may be extended as solutions of equation \eqref{GrindEQ__1_} on \textit{$x<0$} and\textit{$x\ge 0$, }respectively. Then we get

\[f_{2}^{+} \left(x,\lambda \right)=C_{11} \left(\lambda \right)f_{1}^{+} \left(x,\lambda \right)+C_{12} \left(\lambda \right)f_{1}^{-} \left(x,\lambda \right)\, \, \, for\, \, \, x\ge 0,\]

\[f_{1}^{+} \left(x,\lambda \right)=C_{22} \left(\lambda \right)f_{2}^{+} \left(x,\lambda \right)+C_{21} \left(\lambda \right)f_{2}^{-} \left(x,\lambda \right)\, ,\, \, \, \, \, \, \, for\, \, \, x<0\, \, ,\]

\noindent where

\noindent

\[f_{1,2}^{-} \left(x,\lambda \right)=f_{1,2}^{+} \left(x,-\lambda \right),\]

\begin{equation} \label{GrindEQ__10_}
C_{11} \left(\lambda \right)=\frac{W[f_{2}^{+} \left(0,\lambda \right),f_{1}^{-} \left(0,\lambda \right)]}{2i\lambda } ,
\end{equation}

\[C_{12} \left(\lambda \right)=\frac{W[f_{1}^{+} \left(0,\lambda \right),f_{2}^{+} \left(0,\lambda \right)]}{2i\lambda } ,\]

\begin{equation} \label{GrindEQ__11_}
C_{22} \left(\lambda \right)=\frac{i}{\beta } C_{11} \left(-\lambda \right),  C_{21} \left(\lambda \right)=-\frac{i}{\beta } C_{12} \left(\lambda \right).
\end{equation}
\textbf{Proof:} It is easy to see that equation \eqref{GrindEQ__1_} has fundamental solutions$f_{1}^{+} \left(x,\lambda \right)$,$f_{1}^{-} \left(x,\lambda \right)$ ($f_{2}^{+} \left(x,\lambda \right)$,$f_{2}^{-} \left(x,\lambda \right)$) on the$\left|Im\lambda \right|<\frac{\varepsilon }{2} \, \, \, \, (\left|Re\lambda \right|<\frac{\varepsilon }{2} )$, for which

\[W\left[f_{1}^{+} (x,\lambda ),f_{1}^{-} (x,\lambda )\right]=2i\lambda ,\]

\[W\left[f_{2}^{+} (x,\lambda ),f_{2}^{-} (x,\lambda )\right]=2\lambda \beta ,\]
 is satisfied\textbf{}

\noindent  Really, since $W[{\rm \; }f_{1}^{+} \left(x,\lambda \right),f_{1}^{-} \left(x,\lambda \right)]$ and $W[{\rm \; }f_{2}^{+} \left(x,\lambda \right),f_{2}^{-} \left(x,\lambda \right)]$ are independent of $x$ and the functions $f_{1}^{+} \left(x,\lambda \right)$,$f_{1}^{-} \left(x,\lambda \right)$ and $f_{2}^{+} \left(x,\lambda \right)$,$f_{2}^{-} \left(x,\lambda \right)$ allow holomorphic continuation on $x$ to upper and lower half-planes, respectively, the Wronskian coincides as $Imx\to \infty $. We can show that

\begin{equation} \label{GrindEQ__12_}
\mathop{\lim }\limits_{Imx\to \infty } f_{1}^{\pm \left(j\right)} \left(x,\lambda \right)e^{\mp i\lambda x} =\left(\pm i\lambda \right)^{j} \, \, \, \, \, \, \, \, j=0,1,
\end{equation}

\begin{equation} \label{GrindEQ__13_}
\mathop{\lim }\limits_{Imx\to \infty } f_{2}^{\pm \left(j\right)} \left(x,\lambda \right)e^{\mp \lambda x} =\left(\pm \lambda \beta \right)^{j} \, \, \, \, \, \, \, \, j=0,1.
\end{equation}
So that

\[W\left[f_{1}^{+} (x,\lambda ),f_{1}^{-} (x,\lambda )\right]=2i\lambda ,\]

\[W\left[f_{2}^{+} (x,\lambda ),f_{2}^{-} (x,\lambda )\right]=2\lambda \beta .\]
Then each solution of equation \eqref{GrindEQ__1_} may be represented as linear combinations of these solutions.

\[f_{2}^{+} \left(x,\lambda \right)=C_{11} \left(\lambda \right)f_{1}^{+} \left(x,\lambda \right)+C_{12} \left(\lambda \right)f_{1}^{-} \left(x,\lambda \right)\, \, \, for\, \, \, x\ge 0.\]

\[f_{1}^{+} \left(x,\lambda \right)=C_{22} \left(\lambda \right)f_{2}^{+} \left(x,\lambda \right)+C_{21} \left(\lambda \right)f_{2}^{-} \left(x,\lambda \right)\, ,\, \, \, \, \, \, \, for\, \, \, x<0\, \, ,\]
Using the conjunction conditions \eqref{GrindEQ__9_} it is easy to obtain the relation (10-11).

\noindent Let

\begin{equation} \label{GrindEQ__14_}
f_{n}^{\pm } (x)=\mathop{\lim }\limits_{\lambda \to \mp \frac{n}{2} } (n\pm 2\lambda )f_{1}^{\pm } (x,\lambda )=\sum _{\alpha =n}^{\infty }V_{n\alpha } e^{i\alpha x} e^{-i\frac{n}{2} x}  ,
\end{equation}
    It follows from relation \eqref{GrindEQ__6_} that if$V_{nn} \ne 0$, then $V_{n\alpha } \ne 0$ for all $\alpha >n$ and therefore$f_{n}^{\pm } (x)\ne 0$. Consequently, the points $\pm \frac{n}{2} ,\, n\in N$ are not singular points for$f_{1}^{\pm } \left(x,\lambda \right)$.

\noindent Then $W[f_{n}^{\pm } \left(x\right),f_{1}^{\mp } \left(x,\mp \frac{n}{2} \right)]=0$ and consequently the functions$f_{n}^{\pm } \left(x\right),f_{1}^{\mp } \left(x,\mp \frac{n}{2} \right)$, that are solutions of equation \eqref{GrindEQ__1_} for$\lambda =\pm \frac{n}{2} \, \, \, $, are linear dependent.

\noindent Therefore

\begin{equation} \label{GrindEQ__15_}
f_{n}^{\pm } \left(x\right)=V_{nn} f_{1}^{\mp } \left(x,\mp \frac{n}{2} \right),
\end{equation}

\noindent \textbf{\underbar{2.1.  SPECTRUM OF OPERATOR }}$L$.\textbf{\underbar{}}

Let $L$ be an operator generated by the operation  $\frac{1}{\rho \left(x\right)} \left\{-\frac{d^{2} }{dx^{2} } +q\left(x\right)\right\}$

\noindent in the space $L_{2} \left(-\infty ,+\infty ,\rho \left(x\right)\right)$.\textbf{}

\noindent To study the spectrums of the operator $L$ at first we calculate the kernel of the resolvent of the operator  \textbf{$\left(L-\lambda ^{2} I\right)$ }by means of general methods.

\noindent To construct the kernel of the resolvent of operator$L$, we consider the equation

\[-y''\left(x\right)+\, q\left(x\right)y\left(x\right)=\lambda ^{2} \rho \left(x\right)y\left(x\right)+f\left(x\right).\]

Here, $f\left(x\right)$ is an arbitrary function belonging to$L_{2} \left(-\infty ,+\infty \right)$.  Divide the plane $\lambda $ into sectors

\[S_{k} =\{ {\raise0.7ex\hbox{$ k\pi  $}\!\mathord{\left/{\vphantom{k\pi  2}}\right.\kern-\nulldelimiterspace}\!\lower0.7ex\hbox{$ 2 $}} <\arg \lambda <{\raise0.7ex\hbox{$ (k+1)\pi  $}\!\mathord{\left/{\vphantom{(k+1)\pi  2}}\right.\kern-\nulldelimiterspace}\!\lower0.7ex\hbox{$ 2 $}} \} ,k=\overline{0,3}.\]
When $\lambda \in S_{0} $ , we note that every solution of equation \eqref{GrindEQ__1_} can be written in the form

\begin{equation} \label{GrindEQ__16_}
y\left(x,\lambda \right)=C_{1} \left(x,\lambda \right)f_{1}^{+} \left(x,\lambda \right)+C_{2} \left(x,\lambda \right)f_{2}^{+} \left(x,\lambda \right)\, \, \, .
\end{equation}
Using the method of variation of constant, we obtain that

\[C_{1}^{/} \left(x,\lambda \right)=-\frac{1}{W[f_{1}^{+} ,f_{2}^{+} ]} \rho \left(x\right)f_{2}^{+} \left(x,\lambda \right)f\left(x\right)\]

\[C_{2}^{/} \left(x,\lambda \right)=\frac{1}{W[f_{1}^{+} ,f_{2}^{+} ]} \rho \left(x\right)f_{1}^{+} \left(x,\lambda \right)f\left(x\right)\]
By virtue of the condition$y\left(x,\lambda \right)\in L_{2} \left(-\infty ,+\infty \right)$, we find that

\[C_{2} \left(\infty ,\lambda \right)=C_{1} \left(-\infty ,\lambda \right)=0.     \]
Consequently, we have

\[C_{1} \left(x,\lambda \right)=\int _{-\infty }^{x}\frac{1}{W[f_{1}^{+} ,f_{2}^{+} ]} \rho \left(t\right)f_{2}^{+} \left(t,\lambda \right)f\left(t\right) dt\]

\[C_{2} \left(x,\lambda \right)=-\int _{x}^{\infty }\frac{1}{W[f_{1}^{+} ,f_{2}^{+} ]} \rho \left(t\right)f_{1}^{+} \left(t,\lambda \right)f\left(t\right) dt.\]
Substitute them in \eqref{GrindEQ__16_} we get

\[y\left(x,\lambda \right)=\int _{-\infty }^{\infty }R_{11} \left(x,t,\lambda \right) \rho \left(t\right)f\left(t\right)dt\]
where

\begin{equation} \label{GrindEQ__17_}
R_{11} \left(x,t,\lambda \right)=\frac{1}{W[f_{1}^{+} ,f_{2}^{+} ]} \left\{\begin{array}{c} {f_{1}^{+} \left(x,\lambda \right)f_{2}^{+} \left(t,\lambda \right)\, \, \, \, \, \, \, \, \, for\, \, t<x} \\ {f_{1}^{+} \left(t,\lambda \right)f_{2}^{+} \left(x,\lambda \right)\, \, \, \, \, \, \, \, \, for\, \, t>x} \end{array}\right.        \lambda \in S_{0} .
\end{equation}
Calculating analogously we can\textbf{ }construct the kernel of the resolvent on the sectors$S_{k} ,\, \, \, k=\overline{1,3}$, namely

\begin{equation} \label{GrindEQ__18_}
R_{12} \left(x,t,\lambda \right)=\frac{1}{W[f_{1}^{+} ,f_{2}^{-} ]} \left\{\begin{array}{c} {f_{1}^{+} \left(x,\lambda \right)f_{2}^{-} \left(t,\lambda \right)\, \, \, \, \, \, \, \, \, for\, \, t<x} \\ {f_{1}^{+} \left(t,\lambda \right)f_{2}^{-} \left(x,\lambda \right)\, \, \, \, \, \, \, \, \, for\, \, t>x} \end{array}\right.        \lambda \in S_{1} .
\end{equation}

\begin{equation} \label{GrindEQ__19_}
R_{21} \left(x,t,\lambda \right)=\frac{1}{W[f_{1}^{-} ,f_{2}^{-} ]} \left\{\begin{array}{c} {f_{1}^{-} \left(x,\lambda \right)f_{2}^{-} \left(t,\lambda \right)\, \, \, \, \, \, \, \, \, for\, \, t<x} \\ {f_{1}^{-} \left(t,\lambda \right)f_{2}^{-} \left(x,\lambda \right)\, \, \, \, \, \, \, \, \, for\, \, t>x} \end{array}\right.        \lambda \in S_{2} .
\end{equation}

\begin{equation} \label{GrindEQ__20_}
R_{22} \left(x,t,\lambda \right)=\frac{1}{W[f_{1}^{-} ,f_{2}^{+} ]} \left\{\begin{array}{c} {f_{1}^{-} \left(x,\lambda \right)f_{2}^{+} \left(t,\lambda \right)\, \, \, \, \, \, \, \, \, for\, t<x} \\ {f_{1}^{-} \left(t,\lambda \right)f_{2}^{+} \left(x,\lambda \right)\, \, \, \, \, \, \, \, for\, t>x} \end{array}\right.        \lambda \in S_{3} .
\end{equation}

\noindent \textbf{ \underbar{Lemma 2.}   }$L$ has no eigenvalues for real  and pure imaginary $\lambda $. It's continuous

\noindent spectra consist of axes   $Re\lambda =0$  and  $Im\lambda =0$ on which there may exist spectral singularities coinciding with the numbers $\frac{in}{2\beta } ,\, \, \, \frac{n}{2} ,\, \, \, n=\pm 1,\pm 2,\pm 3,...$

\noindent \textbf{Proof:  }We recall that equation \eqref{GrindEQ__1_}\textbf{ }has fundamental solutions$f_{1}^{+} \left(x,\lambda \right)$,$f_{1}^{-} \left(x,\lambda \right)$ ($f_{2}^{+} \left(x,\lambda \right)$,$f_{2}^{-} \left(x,\lambda \right)$) on$\left|Im\lambda \right|<\frac{\varepsilon }{2} \, \, \, \, (\left|Re\lambda \right|<\frac{\varepsilon }{2} )$, for which

\[W\left[f_{1}^{+} (x,\lambda ),f_{1}^{-} (x,\lambda )\right]=2i\lambda ,\]

\[W\left[f_{2}^{+} (x,\lambda ),f_{2}^{-} (x,\lambda )\right]=2\lambda \beta , \]
is satisfied\textbf{}

\noindent Then for $Im\lambda =0$  solution of equation \eqref{GrindEQ__1_} can be written in the form

\[y\left(x,\lambda \right)=C_{1} f_{1}^{+} \left(x,\lambda \right)+C_{2} f_{1}^{-} \left(x,\lambda \right)\, \, \, \]
In case $Im\lambda =0$ the solution  $f_{1}^{\pm } \left(x,\lambda \right)$  has the form

\[f_{1}^{\pm } (x,\lambda )=e^{\pm iRe\lambda x} \left(1+\sum _{n=1}^{\infty }\frac{1}{n\pm 2\lambda } \sum _{\alpha =n}^{\infty }V_{n\alpha } e^{i\alpha x}   \right)\, \, \, \, \, \, \, \, \]
then  $y\left(x,\lambda \right)\in L_{2} \left(-\infty ,+\infty \right)$  except when $C_{1} =C_{2} =0$.

\noindent Analogously we can prove case$Re\lambda =0$. Since in $\left|Re\lambda \right|<\frac{\varepsilon }{2} $  the functions$f_{2}^{+} \left(x,\lambda \right)$, $f_{2}^{-} \left(x,\lambda \right)$ form fundamentals solutions, then

\[y\left(x,\lambda \right)=C_{3} f_{2}^{+} \left(x,\lambda \right)+C_{4} f_{2}^{-} \left(x,\lambda \right)\, \, \, .\]
If  $Re\lambda =0$ then the solution $f_{2}^{\pm } \left(x,\lambda \right)$ has the form

\[f_{2}^{\pm } (x,\lambda )=e^{\pm iIm\lambda \beta x} \left(1+\sum _{n=1}^{\infty }\frac{1}{n\pm 2Im\lambda \beta } \sum _{\alpha =n}^{\infty }V_{n\alpha } e^{i\alpha x}   \right)\, \]
then  $y\left(x,\lambda \right)\in L_{2} \left(-\infty ,+\infty \right)$  except when $C_{3} =C_{4} =0$.

\noindent In order all numbers from the axes $\left\{\lambda :Re\lambda =0\right\}$ and $\left\{\lambda :Im\lambda =0\right\}$ belong to the continuous spectra it suffices to show that domain of value of the operator $\left(L-\lambda ^{2} I\right)$ is dense in $L_{2} \left(-\infty ,+\infty \right)$,  so that the orthogonal complement  of the set $R\left(x,t,\lambda \right)$ consists of only zero element .

\noindent  Let  $\psi \left(x\right)\in L_{2} \left(-\infty ,+\infty \right)$ , $\psi \left(x\right)\ne 0$  and

\begin{equation} \label{GrindEQ__21_}
\int _{-\infty }^{+\infty }\left(Lf-\lambda ^{2} f\right)\overline{\psi \left(x\right)} dx=0
\end{equation}

\noindent be satisfied for any$f\left(x\right)\in D\left(L\right)$.

\noindent From \eqref{GrindEQ__21_} it follows that  $\psi \left(x\right)\in D\left(L^{*} \right)$ and $\psi \left(x\right)$ is an eigenfunction of operator $L^{*} $ corresponding to eigenvalues $\lambda $. More exactly $\overline{\psi \left(x\right)}$  is the solution of the equation

\begin{equation} \label{GrindEQ__22_}
-z''+q\left(x\right)z=\lambda ^{2} z
\end{equation}
belonging to$L_{2} \left(-\infty ,+\infty \right)$.  We obtained that$\psi \left(x\right)=0$, since the operator generated by expression standing at the left hand of \eqref{GrindEQ__22_}, is an operator of type$L$.  This contradiction shows that domain of value of the operator $\left(L-\lambda ^{2} I\right)$ everywhere dense in$L_{2} \left(-\infty ,+\infty \right)$.

\noindent \textbf{Lemma 3.}   The coefficient  $C_{12} \left(\lambda \right)$  is an analytical function in the sector$S_{0} $ and there has finite number of zeros

\noindent \textbf{Proof. }  For solutions  $f_{1}^{\pm } \left(x,\lambda \right)$ and $f_{2}^{\pm } \left(x,\lambda \right)$  we can obtain the asymptotic equalities

\noindent $ $$f_{1}^{\pm \left(j\right)} \left(0,\lambda \right)=\pm \left(i\lambda \right)^{j} +o\eqref{GrindEQ__1_}$         for$\left|\lambda \right|\to \infty ,\, \, \, j=0,1$,

\noindent              for$\left|\lambda \right|\to \infty ,\, \, \, j=0,1$.

\noindent  For simplicity we prove the first equality.

\noindent Since

\[f_{1}^{\pm } \left(0,\lambda \right)=1+\sum _{n=1}^{\infty }\sum _{\alpha =n}^{\infty }\frac{V_{n\alpha } }{n\pm 2\lambda }   \]
that

\[\left|f_{1}^{\pm } \left(0,\lambda \right)\right|\le 1+\sum _{n=1}^{\infty }\sum _{\alpha =n}^{\infty }\frac{\left|V_{n\alpha } \right|}{\left|n+2\lambda \right|} \le   \, 1+\sum _{n=1}^{\infty }\sum _{\alpha =n}^{\infty }\frac{\left|V_{n\alpha } \right|}{\sqrt{\left(n+2Re\lambda \right)^{2} +4Im^{2} \lambda } } \le   \, \, 1+\frac{1}{\left|Im\lambda \right|} \sum _{n=1}^{\infty }\sum _{\alpha =n}^{\infty }\frac{\alpha \left|V_{n\alpha } \right|}{n}   .\]
Therefore, as$\left|\lambda \right|\to \infty $, we obtain$f_{1}^{\pm } \left(0,\lambda \right)=1+o\left(1\right)$.

\noindent Analogously we can prove the rest asymptotic equalities as$\left|\lambda \right|\to \infty $,

\noindent Then for the coefficients $C_{12} (\lambda ),\, \, C_{11} (-\lambda ),\, \, C_{12} (-\lambda ),\, \, C_{11} (\lambda )$ we get the following asymptotic equalities

\[C_{12} \left(\lambda \right)=\frac{1}{2i\lambda } \left(\lambda \beta -i\lambda \right)+o\left(1\right)=-\frac{i\beta +1}{2} +o\left(1\right),\]

\[C_{12} \left(-\lambda \right)=-\frac{i\beta +1}{2} +o\left(1\right),\]

\[C_{11} \left(\lambda \right)=-\frac{1-i\beta }{2} +o\left(1\right),\]

\[C_{11} \left(-\lambda \right)=-\frac{1-i\beta }{2} +o\left(1\right).\]
These asymptotic equalities and analytical properties of the coefficients $C_{12} (\lambda ),\, \, C_{11} (-\lambda ),\, \, C_{12} (-\lambda ),\, \, C_{11} (\lambda )$  make valid the following statement.

\noindent \textbf{Lemma 4.}  The eigenvalues of operator $L$are finite and coincide with zeros of the functions $C_{12} (\lambda ),\, \, C_{11} (-\lambda ),\, \, C_{12} (-\lambda ),\, \, C_{11} (\lambda )$ from sectors

\[S_{k} =\{ {\raise0.7ex\hbox{$ k\pi  $}\!\mathord{\left/{\vphantom{k\pi  2}}\right.\kern-\nulldelimiterspace}\!\lower0.7ex\hbox{$ 2 $}} <\arg \lambda <{\raise0.7ex\hbox{$ (k+1)\pi  $}\!\mathord{\left/{\vphantom{(k+1)\pi  2}}\right.\kern-\nulldelimiterspace}\!\lower0.7ex\hbox{$ 2 $}} \} ,k=\overline{0,3}\]
respectively.

\noindent \textbf{Definition 1.}  The data $\{ \, \lambda _{n} ,\, C_{11} \left(\lambda \right),\, \, C_{12} \left(\lambda \right)\} $are called the spectral data of$L$.

\noindent

\noindent \textbf{2.2. EIGENFUNCTION EXPANSIONS.}

\noindent \textbf{}

\noindent \textbf{Definition 2}.  The points at which resolvent have poles are called the singular numbers of operator$L$.

\noindent Let  $\lambda _{1} ,\lambda _{2} ,....\lambda _{l} ,\lambda _{l+1} .....\lambda _{n} ...$be the singular numbers of operator $L$.At that

\[Re\lambda _{j} Im\lambda _{j} \ne 0,\, \, \, \, \, \, \, j=1,2,....l\]

\[Re\lambda _{j} Im\lambda _{j} =0,\, \, \, \, \, \, \, j=l+1,....n,...\]
 Order$k_{j} $ of root$\lambda _{j} $is called the order of root of the singular numbers$\lambda _{j} ,\, \, \, j=1,..l.$.  It is clear that then the numbers $\lambda _{1,} \lambda _{2,} ....\lambda _{l,} $ will be eigenvalues of operator$L$.The numbers $\lambda _{j} ,\, \, \, j=l+1,....n,..$are called the spectral singularities of operator$L$. From the form of resolvent it is easy to see that it has singular numbers (i.e. eigenvalues) $\lambda _{1,} \lambda _{2,} ....\lambda _{l,} $in zeros of the functions  $C_{12} (\lambda ),\, \, C_{11} (-\lambda ),\, \, C_{12} (-\lambda ),\, \, C_{11} (\lambda )$ in the sectors $S_{k} =\{ {\raise0.7ex\hbox{$ k\pi  $}\!\mathord{\left/{\vphantom{k\pi  2}}\right.\kern-\nulldelimiterspace}\!\lower0.7ex\hbox{$ 2 $}} <\arg \lambda <{\raise0.7ex\hbox{$ (k+1)\pi  $}\!\mathord{\left/{\vphantom{(k+1)\pi  2}}\right.\kern-\nulldelimiterspace}\!\lower0.7ex\hbox{$ 2 $}} \} ,k=\overline{0,3}$  respectively. Their finiteness follows from Lemma 4. It directly follows from Lemma 2 and representation (17-20) that kernel of resolvent may have spectral singularities coinciding with the numbers$\frac{in}{2\beta } ,\, \, \, \frac{n}{2} ,\, \, \, n=\pm 1,\pm 2,\pm 3,...$.  Consequently taking \eqref{GrindEQ__15_} into account and using

\[\mathop{\lim }\limits_{\lambda \to {\raise0.7ex\hbox{$ n $}\!\mathord{\left/{\vphantom{n 2}}\right.\kern-\nulldelimiterspace}\!\lower0.7ex\hbox{$ 2 $}} } \left(n-2\lambda \right)W[f_{2}^{+} \left(0,\lambda \right),f_{1}^{-} \left(0,\lambda \right)]=V_{nn} W[f_{2}^{+} \left(0,\frac{n}{2} \right),f_{1}^{+} \left(0,\frac{n}{2} \right)]\]
we calculate

\begin{equation} \label{GrindEQ__23_}
\begin{array}{l} {\mathop{\lim }\limits_{\lambda \to {\raise0.7ex\hbox{$ n $}\!\mathord{\left/{\vphantom{n 2}}\right.\kern-\nulldelimiterspace}\!\lower0.7ex\hbox{$ 2 $}} } \left(n-2\lambda \right)R_{11} \left(x,t,\lambda \right)=\mathop{\lim }\limits_{\lambda \to {\raise0.7ex\hbox{$ n $}\!\mathord{\left/{\vphantom{n 2}}\right.\kern-\nulldelimiterspace}\!\lower0.7ex\hbox{$ 2 $}} } \left(n-2\lambda \right)\frac{1}{2i\lambda } [f_{1}^{+} \left(x,\lambda \right)f_{1}^{+} \left(t,\lambda \right)\frac{W[f_{2}^{+} ,f_{1}^{-} ]}{W[f_{1}^{+} ,f_{2}^{+} ]} +} \\ {+f_{1}^{+} \left(x,\lambda \right)f_{1}^{-} \left(t,\lambda \right)]=\frac{1}{in} [V_{nn} f_{1}^{+} \left(x,\frac{n}{2} \right)f_{1}^{+} \left(t,\frac{n}{2} \right)+} \\ {+V_{nn} f_{1}^{+} \left(x,\frac{n}{2} \right)f_{1}^{+} \left(t,\frac{n}{2} \right)]=\frac{2}{in} V_{nn} f_{1}^{+} \left(x,\frac{n}{2} \right)f_{1}^{+} \left(t,\frac{n}{2} \right).} \end{array}
\end{equation}
Analogously taking into account

\[\mathop{\lim }\limits_{\lambda \to {\raise0.7ex\hbox{$ in $}\!\mathord{\left/{\vphantom{in 2\beta }}\right.\kern-\nulldelimiterspace}\!\lower0.7ex\hbox{$ 2\beta  $}} } \left(n+2i\lambda \beta \right)f_{2}^{-} \left(x,\lambda \right)=V_{nn} f_{2}^{+} \left(x,\frac{in}{2\beta } \right)\]

\[\mathop{\lim }\limits_{\lambda \to {\raise0.7ex\hbox{$ in $}\!\mathord{\left/{\vphantom{in 2\beta }}\right.\kern-\nulldelimiterspace}\!\lower0.7ex\hbox{$ 2\beta  $}} } \left(n+2i\lambda \beta \right)W[f_{2}^{-} \left(0,\lambda \right),f_{1}^{+} \left(0,\lambda \right)]=V_{nn} W[f_{2}^{+} \left(0,\frac{in}{2\beta } \right),f_{1}^{+} \left(0,\frac{in}{2\beta } \right)]\]
we get

\begin{equation} \label{GrindEQ__24_}
\begin{array}{l} {\mathop{\lim }\limits_{\lambda \to {\raise0.7ex\hbox{$ in $}\!\mathord{\left/{\vphantom{in 2}}\right.\kern-\nulldelimiterspace}\!\lower0.7ex\hbox{$ 2 $}} \beta } \left(n+2i\lambda \beta \right)R_{11} \left(x,t,\lambda \right)=\mathop{\lim }\limits_{\lambda \to {\raise0.7ex\hbox{$ in $}\!\mathord{\left/{\vphantom{in 2}}\right.\kern-\nulldelimiterspace}\!\lower0.7ex\hbox{$ 2 $}} \beta } \left(n+2i\lambda \beta \right)\frac{1}{2\lambda \beta } [f_{2}^{+} \left(x,\lambda \right)f_{2}^{+} \left(t,\lambda \right)\frac{W[f_{2}^{-} ,f_{1}^{+} ]}{W[f_{2}^{+} ,f_{1}^{+} ]} +} \\ {+f_{2}^{+} \left(x,\lambda \right)f_{2}^{-} \left(t,\lambda \right)]==\frac{1}{in} [V_{nn} f_{2}^{+} \left(x,\frac{in}{2\beta } \right)f_{2}^{+} \left(t,\frac{in}{2\beta } \right)+} \\ {+V_{nn} f_{2}^{+} \left(x,\frac{in}{2\beta } \right)f_{2}^{+} \left(t,\frac{in}{2\beta } \right)]=\frac{2}{in} V_{nn} f_{2}^{+} \left(x,\frac{in}{2\beta } \right)f_{2}^{+} \left(t,\frac{in}{2\beta } \right).} \end{array}
\end{equation}

\noindent

\noindent

\noindent

\noindent

\noindent Calculate rest residues.

\begin{equation} \label{GrindEQ__25_}
\begin{array}{l} {\mathop{\lim }\limits_{\lambda \to -{\raise0.7ex\hbox{$ n $}\!\mathord{\left/{\vphantom{n 2}}\right.\kern-\nulldelimiterspace}\!\lower0.7ex\hbox{$ 2 $}} } \left(n+2\lambda \right)R_{12} \left(x,t,\lambda \right)=\mathop{\lim }\limits_{\lambda \to -{\raise0.7ex\hbox{$ n $}\!\mathord{\left/{\vphantom{n 2}}\right.\kern-\nulldelimiterspace}\!\lower0.7ex\hbox{$ 2 $}} } \left(n+2\lambda \right)\frac{1}{2i\lambda } [f_{1}^{+} \left(x,\lambda \right)f_{1}^{+} \left(t,\lambda \right)\frac{W[f_{1}^{-} ,f_{2}^{-} ]}{W[f_{2}^{-} ,f_{1}^{+} ]} +} \\ {+f_{1}^{+} \left(x,\lambda \right)f_{1}^{-} \left(t,\lambda \right)]=\frac{1}{in} [V_{nn} f_{1}^{-} \left(x,-\frac{n}{2} \right)\, \, \tilde{f}_{1}^{+} \left(t,\frac{n}{2} \right)\frac{W[f_{1}^{-} \left(0,-\frac{n}{2} \right),f_{2}^{-} \left(0,-\frac{n}{2} \right)]}{W[f_{2}^{-} \left(0,-\frac{n}{2} \right),\tilde{f}_{1}^{+} \left(0,-\frac{n}{2} \right)]} +} \\ {+V_{nn} f_{1}^{-} \left(x,-\frac{n}{2} \right)f_{1}^{-} \left(t,-\frac{n}{2} \right)]=V_{nn} F_{1} \left(x,t,-\frac{n}{2} \right).} \end{array}
\end{equation}
Here we use denotation$\tilde{f}_{1}^{+} \left(x,\lambda \right)=f_{1}^{+} \left(x,\lambda \right)\left(n+2\lambda \right)$, therewith, the function $\tilde{f}_{1}^{+} \left(x,\lambda \right)$ has no poles at the points$\lambda =-\frac{n}{2} ,\, \, \, n\in N$, and  $F_{1} \left(x,t,-\frac{n}{2} \right)=\frac{1}{in} [f_{1}^{-} \left(x,-\frac{n}{2} \right)\, \, \tilde{f}_{1}^{+} \left(t,\frac{n}{2} \right)\frac{W[f_{1}^{-} \left(0,-\frac{n}{2} \right),f_{2}^{-} \left(0,-\frac{n}{2} \right)]}{W[f_{2}^{-} \left(0,-\frac{n}{2} \right),\tilde{f}_{1}^{+} \left(0,-\frac{n}{2} \right)]} +f_{1}^{-} \left(x,-\frac{n}{2} \right)f_{1}^{-} \left(t,-\frac{n}{2} \right)].$

\noindent Analogously we denote$\tilde{f}_{2}^{+} \left(x,\lambda \right)=f_{2}^{+} \left(x,\lambda \right)\left(n-2i\lambda \beta \right)$. Remark that the function $\tilde{f}_{2}^{+} \left(x,\lambda \right)$ also has no poles at the points$\lambda =-\frac{in}{2\beta } ,\, \, \, n\in N$, and

\[\begin{array}{l} {F_{2} \left(x,t,-\frac{in}{2\beta } \right)=\frac{1}{in} [f_{2}^{-} \left(x,-\frac{in}{2\beta } \right)f_{2}^{-} \left(t,-\frac{in}{2\beta } \right)+} \\ {+\frac{W[f_{2}^{-} \left(0,-\frac{in}{2\beta } \right),f_{1}^{-} \left(0,-\frac{in}{2\beta } \right)]}{W[\tilde{f}_{2}^{+} \left(0,-\frac{in}{2\beta } \right),f_{1}^{-} \left(0,-\frac{in}{2\beta } \right)]} \tilde{f}_{2}^{+} \left(0,-\frac{in}{2\beta } \right)f_{2}^{-} \left(t,-\frac{in}{2\beta } \right)].} \end{array}\]
Then we get

\[\begin{array}{l} {\mathop{\lim }\limits_{\lambda \to -{\raise0.7ex\hbox{$ in $}\!\mathord{\left/{\vphantom{in 2}}\right.\kern-\nulldelimiterspace}\!\lower0.7ex\hbox{$ 2 $}} \beta } \left(n-2i\lambda \beta \right)R_{22} \left(x,t,\lambda \right)=-\mathop{\lim }\limits_{\lambda \to -{\raise0.7ex\hbox{$ in $}\!\mathord{\left/{\vphantom{in 2\beta }}\right.\kern-\nulldelimiterspace}\!\lower0.7ex\hbox{$ 2\beta  $}} } \left(n-2i\lambda \beta \right)\frac{1}{2\lambda \beta } [f_{2}^{-} \left(x,\lambda \right)f_{2}^{+} \left(t,\lambda \right)\frac{W[f_{2}^{-} ,f_{1}^{-} ]}{W[f_{2}^{+} ,f_{1}^{-} ]} +} \\ {+f_{2}^{+} \left(x,\lambda \right)f_{2}^{+} \left(t,\lambda \right)]=\frac{1}{in} [V_{nn} f_{2}^{-} \left(x,-\frac{in}{2\beta } \right)f_{2}^{-} \left(t,-\frac{in}{2\beta } \right)+} \\ {+V_{nn} \frac{W[f_{2}^{-} \left(0,-\frac{in}{2\beta } \right),f_{1}^{-} \left(0,-\frac{in}{2\beta } \right)]}{W[\tilde{f}_{2}^{+} \left(0,-\frac{in}{2\beta } \right),f_{1}^{-} \left(0,-\frac{in}{2\beta } \right)]} \tilde{f}_{2}^{+} \left(0,-\frac{in}{2\beta } \right)f_{2}^{-} \left(t,-\frac{in}{2\beta } \right)]=V_{nn} F_{2} \left(x,t,-\frac{in}{2\beta } \right).} \end{array}\]
\textbf{Lemma 6:}   Let  $\psi \left(x\right)$ be an arbitrary twice continuously differentiable function belonging to$L_{2} \left(-\infty ,+\infty ,\rho \left(x\right)\right)$ . Then

\[\int _{-\infty }^{+\infty }R\left(x,t,\lambda \right) \rho \left(t\right)\psi \left(t\right)dt=-\frac{\psi \left(x\right)}{\lambda ^{2} } +\frac{1}{\lambda ^{2} } \int _{-\infty }^{+\infty }R\left(x,t,\lambda \right) g\left(t\right)dt,\]
  where

\[g\left(t\right)=-\psi ''\left(x\right)+q\left(x\right)\psi \left(x\right)\in L_{2} \left(-\infty ,+\infty \right).\]
Integrating the both hand side along the circle $\left|\lambda \right|=R$ and passing to limit as $R\to \infty $   we get

\[\psi \left(x\right)=-\mathop{\lim }\limits_{R\to \infty } \frac{1}{2\pi i} \oint _{\left|\lambda \right|=R}2\lambda d\lambda \int _{-\infty }^{+\infty }R\left(x,t,\lambda \right)\rho \left(t\right)\psi \left(t\right)dt  \]
The function  $\int _{-\infty }^{+\infty }R\left(x,t,\lambda \right)\rho \left(t\right)\psi \left(t\right)dt $ is analytical inside the contour, with respect to $\lambda $ excepting the points$\lambda =\lambda _{n} ,\, \, \lambda =\pm \frac{n}{2} ,\, \, \lambda =\pm \frac{in}{2\beta } ,\, \, \, n=1,2,...$.   Denote by  $_{0}^{+} \left(_{0}^{-} \right)$   the contour formed by segments  $[0,\frac{1}{2} +\delta ],[\frac{n}{2} +\delta ,\frac{\left(n+1\right)}{2} -\delta ]$   and semi-circles of radius $\delta $ with the centers at points  $\frac{n}{2} ,\, \, \, \, n=1,2,...$ located in upper (lower) half plane. Analogously we denote by  the $_{0i}^{+} \left(_{0i}^{-} \right)$ the contour formed by segments $[0,\frac{i}{2\beta } +\delta ],[\frac{in}{2\beta } +\delta ,\frac{i\left(n+1\right)}{2\beta } -\delta ],\, \, \, \, n=1,2,3....$ and semi-circles of radius $\delta $ with the centers at points  $\, \frac{in}{2\beta } ,\, \, \, n=1,2,...$ located in right (left) half plane.

\noindent  Let the contours  $ \Gamma _1^ +  \left( {\Gamma _1^ -}
\right)$
 and $\Gamma _{1i}^ +  \left( {\Gamma _{1i}^ -  } \right)  $ be
obtained from $ \Gamma _0^ +  \left( {\Gamma _0^ -  } \right) $
and $ \Gamma _{0i}^ +  \left( {\Gamma _{0i}^ -  } \right) $   by
turning around the angle$\pi $. Then

\[\begin{array}{l} {\psi \left(x\right)=-\frac{1}{2i\pi } \int _{-\infty }^{+\infty }2\lambda \rho \left(t\right)\psi \left(t\right) [\int _{\Gamma_{0i}^{+} }R_{11} \left(x,t,\lambda \right) d\lambda +\int _{\Gamma_{0}^{+} }R_{11} \left(x,t,\lambda \right) d\lambda -\int _{\Gamma_{0i}^{-} }R_{12} \left(x,t,\lambda \right) d\lambda +} \\ {+\int _{\Gamma_{1}^{+} }R_{12} \left(x,t,\lambda \right) d\lambda -\int _{\Gamma_{1}^{-} }R_{21} \left(x,t,\lambda \right) d\lambda -\int _{\Gamma_{1i}^{-} }R_{21} \left(x,t,\lambda \right) d\lambda +\int _{\Gamma_{1}^{+} }R_{22} \left(x,t,\lambda \right) d\lambda -} \\ {-\int _{\Gamma_{0}^{-} }R_{21} \left(x,t,\lambda \right) d\lambda ]dt=-\frac{1}{2i\pi } \int _{-\infty }^{+\infty }2\lambda \rho \left(t\right)\psi \left(t\right) [\int _{\Gamma_{0i}^{-} }[R_{11} \left(x,t,\lambda \right) -R_{12} \left(x,t,\lambda \right)]d\lambda +} \\ {+\int _{\Gamma_{0}^{-} }[R_{11} \left(x,t,\lambda \right) -R_{22} \left(x,t,\lambda \right)]d\lambda +\int _{\Gamma_{1}^{-} }[R_{12} \left(x,t,\lambda \right) -R_{21} \left(x,t,\lambda \right)]d\lambda +} \\ {+\int _{\Gamma_{1i}^{-} }[R_{22} \left(x,t,\lambda \right) -R_{21} \left(x,t,\lambda \right)]d\lambda +\mathop{+Res}\limits_{\lambda =\lambda _{n} } [R_{11} \left(x,t,\lambda \right)+R_{12} \left(x,t,\lambda \right)+R_{21} \left(x,t,\lambda \right)+} \\ {+R_{22} \left(x,t,\lambda \right)]\mathop{+Res}\limits_{\lambda =\frac{in}{2\beta } } R_{11} \left(x,t,\lambda \right)\mathop{+Res}\limits_{\lambda =\frac{n}{2} } R_{11} \left(x,t,\lambda \right)\mathop{+Res}\limits_{\lambda =-\frac{n}{2} } R_{12} \left(x,t,\lambda \right)+\mathop{+Res}\limits_{\lambda =-\frac{in}{2\beta } } R_{22} \left(x,t,\lambda \right)]dt} \end{array}\]
Separately calculate every item.

\[R_{11} \left(x,t,\lambda \right)-R_{12} \left(x,t,\lambda \right)=\frac{f_{1}^{+} \left(x,\lambda \right)f_{1}^{+} \left(t,\lambda \right)}{2i\lambda C_{12} \left(\lambda \right)C_{22} \left(\lambda \right)} \]

\[R_{11} \left(x,t,\lambda \right)-R_{22} \left(x,t,\lambda \right)=\frac{f_{2}^{+} \left(x,\lambda \right)f_{2}^{+} \left(t,\lambda \right)}{2i\lambda C_{12} \left(\lambda \right)C_{11} \left(\lambda \right)} \]

\[R_{12} \left(x,t,\lambda \right)-R_{21} \left(x,t,\lambda \right)=\frac{f_{2}^{-} \left(x,\lambda \right)f_{2}^{-} \left(t,\lambda \right)}{2i\lambda C_{11} \left(-\lambda \right)C_{12} \left(-\lambda \right)} \]

\[R_{22} \left(x,t,\lambda \right)-R_{21} \left(x,t,\lambda \right)=\frac{f_{1}^{-} \left(x,\lambda \right)f_{1}^{-} \left(t,\lambda \right)}{2i\lambda C_{22} \left(-\lambda \right)C_{21} \left(-\lambda \right)} \]

\noindent  Residues of resolvent $R_{pq} \, \, \, p,q=1,2$ in $\lambda _{1,} \lambda _{2,} ....\lambda _{l,} $ denote by $G_{pq} \left(\lambda _{n} \right)$ . Then $G_{pq} \left(\lambda _{n} \right)$ will be equal to

\[G_{pq} \left(\lambda _{n} \right)=\frac{1}{\left(k_{n} -1\right)!} \mathop{\lim }\limits_{\lambda \to \lambda _{n} } \frac{d^{k_{n} -1} }{d\lambda ^{k_{n-1} } } [\left(\lambda -\lambda _{n} \right)^{k_{n} } R_{pq} \left(x,t,\lambda \right)]. \]
Then for every function $\psi \left(x\right)$ belonging to $L_{2} \left(-\infty ,+\infty ,\rho \left(x\right)\right)$ we get following eigenfunction expansion in the form

\[\begin{array}{l} {\psi \left(x\right)=-\frac{1}{2i\pi } \int _{-\infty }^{+\infty }\rho \left(t\right)\psi \left(t\right) [\int _{_{0i}^{-} }[\frac{f_{1}^{+} \left(x,\lambda \right)f_{1}^{+} \left(t,\lambda \right)}{iC_{12} \left(\lambda \right)C_{22} \left(\lambda \right)}  ]d\lambda +} \\ {+\int _{_{0}^{-} }[\frac{f_{2}^{+} \left(x,\lambda \right)f_{2}^{+} \left(t,\lambda \right)}{iC_{12} \left(\lambda \right)C_{11} \left(\lambda \right)}  ]d\lambda +\int _{_{1}^{-} }[\frac{f_{2}^{-} \left(x,\lambda \right)f_{2}^{-} \left(t,\lambda \right)}{iC_{11} \left(-\lambda \right)C_{12} \left(-\lambda \right)}  ]d\lambda +} \\ {+\int _{_{1ii}^{-} }[\frac{f_{1}^{-} \left(x,\lambda \right)f_{1}^{-} \left(t,\lambda \right)}{iC_{22} \left(-\lambda \right)C_{21} \left(-\lambda \right)}  ]d\lambda \mathop{+\sum _{p,q=1}^{2}G_{pq}  }\limits_{\lambda =\lambda _{n} } +V_{nn} F\left(x,t,n\right)]dt} \end{array}\]
where

\noindent

\[F\left(x,t,n\right)=\frac{2}{in} f_{1}^{+} \left(x,\frac{n}{2} \right)f_{1}^{+} \left(t,\frac{n}{2} \right)+f_{2}^{+} \left(x,\frac{in}{2\beta } \right)f_{2}^{+} \left(t,\frac{in}{2\beta } \right)+F_{1} \left(x,t,\frac{n}{2} \right)+F_{2} \left(x,t,-\frac{in}{2\beta } \right)\]

\noindent

\noindent \textbf{SOLUTION OF THE INVERSE PROBLEM.}

\noindent \textbf{}

\noindent \textbf{}

\noindent Let's study the inverse problem for the problem (1-3). The inverse problem is formulat as follows.

\noindent \textbf{INVERSE PROBLEM. }Given the spectral data $\{ \, \lambda _{n} ,\, C_{11} \left(\lambda \right),\, \, C_{12} \left(\lambda \right)\} $construct the$\beta $

\noindent and potential $q\left(x\right)$.

\noindent Using the results obtained above we arrive at the following procedure for solution of the inverse problem.

\noindent 1. Taking into account \eqref{GrindEQ__15_} it is easy to check that

\[\mathop{\lim }\limits_{\lambda \to {\raise0.7ex\hbox{$ n $}\!\mathord{\left/{\vphantom{n 2}}\right.\kern-\nulldelimiterspace}\!\lower0.7ex\hbox{$ 2 $}} } \left(n-2\lambda \right)\frac{C_{11} \left(\lambda \right)}{C_{12} \left(\lambda \right)} =V_{nn} ,\]
consequently we find all numbers $V_{nn} ,\, \, n=1,2,...$.

\noindent

\noindent 2. Taking into account \eqref{GrindEQ__6_} we get

\[V_{n,\alpha +n} =V_{nn} \sum _{m=1}^{\alpha }\frac{V_{m\alpha } }{m+n}  ,\]
from which all numbers  $V_{n\alpha } ,\, \, \, \alpha =1,2.....,\, \, \, n=1,2,....n<\alpha $ are defined.

\noindent 3. Then from recurrent formula \eqref{GrindEQ__6_}-\eqref{GrindEQ__8_}, find all numbers$q_{n} $.

\noindent 4. The number $\beta $  is defined from equality

\[\beta =iC_{11} \left(\lambda _{n} \right)C_{11} \left(-\lambda _{n} \right).\]
Really using Lemma 1 we derive for $\lambda _{n} \in S_{0} $ the true relation

\noindent $C_{11} \left(\lambda \right)=\frac{f_{2}^{+} \left(x,\lambda _{n} \right)}{f_{1}^{+} \left(x,\lambda _{n} \right)} $, $C_{22} \left(\lambda \right)=\frac{f_{1}^{+} \left(x,\lambda _{n} \right)}{f_{2}^{+} \left(x,\lambda _{n} \right)} $ i.e. $C_{11} \left(\lambda _{n} \right)C_{22} \left(\lambda _{n} \right)=1$.

\noindent Then from \eqref{GrindEQ__11_} we get  $\beta =iC_{11}
\left(\lambda _{n} \right)C_{11} \left(-\lambda _{n} \right)$.

So inverse problem has a unique solution and the numbers $\beta $ and $q_{n} $ are defined constructively by the spectral data.

\noindent \textbf{Theorem 2.} The specification of the spectral
data uniquely determines $\beta $\,\,\,\,and
potential\,\,$q\left(x\right)$.

\noindent

\noindent \textbf{REFERENCES.}

\noindent 1.\textit{   }Kamke E. Handbook of ordinary Differential equations (Russian).  Nauka,

\noindent       Moscow,1976.

\noindent 2.   Jeffrey C. Lagarias. The Schrödinger Operator with Morse Potential on the

\noindent       right half line. Arxiv: 0712.3238v1 [math.SP] 19 Dec 2007.

\noindent 3    Gasymov M.G. Spectral analysis of a class non-self-adjoint operator of the

\noindent       second order. Functional analysis and its appendix. (In  Russian) 1980,

\noindent       V34.1.pp.14-19

\noindent 4.   Guseinov, I.~M. and Pashaev, R.~T.~ On an inverse problem for a second-

\noindent       order differential equation. UMN, 2002 , 57:3, 147--148

\noindent 5.   V.Guillemin, A. Uribe Hardy functionms and the inverse spectral method.

\noindent       Comm.In Partial Differential equations, 8\eqref{GrindEQ__13_},1455-1474(1983)

\noindent 6.   Pastur L.A., Tkachenko V.A. An inverse problem for one class of

\noindent       onedimentional Shchrodinger's operators with complex periodic potentials.

\noindent       Funksional analysis and its appendix. ( in
Russian) 1990 V54.¹6.pp. 1252-1269
 appendix.(inRussian)1990 V54.¹6.pp.
1252-1269

 \noindent7.    Efendiev, R.F. Spectral analysis of a class of non-self-adjoint
                differential operator pencils with a generalized function.Teoreticheskaya i

\noindent       Matematicheskaya  Fizika, 2005,Vol.145,1.pp.102-107,October(Russian).

\noindent       Theoretical and Mathematical Physics,145\eqref{GrindEQ__1_}:1457-461,(English).

\noindent 8.  Efendiev, R. F. omplete solution of an inverse problem for one class of the

\noindent       high order ordinary differential operators with periodic coefficients. Zh. Mat.

\noindent       Fiz. Anal. Geom2006, C.\textit{ }2, no. 1, 73--86, 111.

\noindent 9. Efendiev, R.F. The Characterization Problem for One Class of Second Order

\noindent       Operator Pencil with Complex Periodic Coefficients. Moscow Mathematical

\noindent       Journal, 2007, Volume~7,~Number~1, pp.55-65.

\end{document}